\documentclass[12pt,leqno,fleqn]{amsart}  
\usepackage{amsmath,amstext,amsthm,amssymb,amsxtra}
\usepackage{mathtools}
\mathtoolsset{showonlyrefs,showmanualtags}

\usepackage{hyperref} 
\hypersetup{
    colorlinks=true,       
    linkcolor=blue,          
    citecolor=magenta,        
    filecolor=magenta,      
    urlcolor=cyan           
}
\usepackage[alphabetic]{amsrefs}

\allowdisplaybreaks

\theoremstyle{plain} 
\newtheorem{lemma}[equation]{Lemma} 
 
\newtheorem{theorem}[equation]{Theorem}

\theoremstyle{definition}
\newtheorem{definition}[equation]{Definition} 
\newtheorem{czo}[equation]{Calder\'on-Zygmund Operators}

\theoremstyle{remark}

\newtheorem*{ack}{Acknowledgment}

\numberwithin{equation}{section}
\numberwithin{paragraph}{section}

\title {The Linear Bound in $ A_2$ for Calder\'on-Zygmund Operators: A Survey}
\author{Michael Lacey}   

\address{ School of Mathematics, Georgia Institute of Technology, Atlanta GA 30332, USA}
\email {lacey@math.gatech.edu}
\thanks{Research supported in part by grant NSF-DMS 0968499.}
\begin{document}

\begin{abstract}
For an $ L ^2 $-bounded Calder\'on-Zygmund Operator $ T$ acting on $ L ^2 (\mathbb R ^{d})$, and a weight $ w \in A_2$, the norm of $ T$ on $ L ^2 (w)$ is dominated by $ C _T \lVert w\rVert_{A_2}$.  The recent theorem completes a line of investigation initiated by Hunt-Muckenhoupt-Wheeden in 1973 \cite{MR0312139},  has been established in different levels of generality by a number of authors over the last few years. 
It has a subtle proof, whose full implications will unfold over the next few years.  
This sharp estimate requires that the $ A_2$ character of the weight can be exactly once in the proof.  
Accordingly, a large part of the proof uses two-weight techniques,  is based on novel decomposition methods for operators and weights, and yields new insights into the Calder\'on-Zygmund theory.  
We survey the proof of this Theorem in this paper.  
\end{abstract}

\maketitle

\section{Introduction} 

We  survey recent developments on the norm behavior of classical Calder\'on-Zygmund operators on weighted spaces, with a special focus on the Muckenhoupt--Wheeden class of weights $ A_2$.   
Indeed, after the 40 some-odd years since the class of $ A_p$ weights was introduced by Muckenhoupt and Wheeden, the theory has reached a natural milestone, with the sharp dependence of norm estimates being established.  
We concentrate on an exposition of the techniques behind this Theorem: 

\begin{theorem}\label{t.linear} Let $ T $ be an $ L ^2 $-bounded Calder\'on-Zygmund operator acting on $ L ^2 (\mathbb R ^{d})$ (for precise definition see Definition~\ref{d.CZ}). 
And, let $ w \in A_2$ (for precise definition see Definition~\ref{d.Ap}).   We then have the estimate 
\begin{equation}\label{e.linear}
\lVert T f \rVert_{ L ^2 (w)} \le C _{T} \lVert w\rVert_{A_2} \lVert f\rVert_{L ^2 (w)} \,. 
\end{equation}
Here $0< C_T < \infty $ depends only on the operator $ T$ and dimension $ d$. 
\end{theorem}

The theory of weights came of age in 1973, with the result of Hunt-Muckenhoupt-Wheeden \cite{MR0312139}, which showed in one dimension that for $ w >0$ a.e., the Hilbert transform is bounded on $ L ^2 (w)$ if and only if $ w\in A_2$.  
This result was established for other suitable collections of singular integrals in higher dimensions.  
But, early proofs of this fact delivered a poor control on the norm, and indeed, the significance of the sharp dependence was a theme recognized over time.  

The interest here is that the power of the $ A_2$ characteristic is in general sharp.  
Accordingly, the method of proof is delicate, and indeed sheds new light on methods and techniques appropriate for weighted spaces, as well as the structure of Calder\'on-Zygmund operators.  

It is known that the estimate \eqref{e.linear}, together with sharp extrapolation \cite{MR2140200}, gives the sharp estimate on $ L ^{p} (w)$, for 
$ 1<p< \infty $, accordingly, we concentrate on the $ L ^2 $ case. 
We recall definitions in the next two sections, and then recall different elements of one of the proofs known of this paper, the pleasingly direct proof of Hyt\"onen-P\'erez-Treil-Volberg \cite{1010.0755}.  
The concluding section includes some historical remarks, and a variety of pointers to cognate results and approaches.  

\begin{ack}
Due to my, and my father's, personal connection to Polish mathematicians, it was my distinct pleasure to participate in the conference 
marking the centenary of birth of J\'ozef Marcinkiewicz. 
It was a fitting testament to the life of Marcinkiewicz, of what was accomplished,  what was lost, and finally, what the people of Poznan and Poland can now achieve, in their beautiful and prosperous city and country.    
\end{ack}

\section{Calder\'on-Zygmund Operators} 

There are two canonical examples of Calder\'on-Zygmund operators that one can keep in mind. 
The first is Hilbert transform itself, defined by 
\begin{equation}\label{e.Hdef}
H f (x) := \lim _{\epsilon \to 0} \int f (x-y) \; \frac {dy} y \,.
\end{equation}
Here, one should note that if $ f$ is Schwartz class, then the limit above exists for all $ x$, and is referred to as the principal value of the integral.  In brief, $ H f (x) = \textup{pv}\; f \ast \frac 1x$.  
But, the Hilbert transform is a convolution operator, which introduces a subtle simplification in its analysis in Lebesgue space.  (There is no paraproduct to control.)
Aside from the Hilbert transform, the other canonical convolution operators are the Beurling in the plane, and the vector of Riesz transforms.

A second example to keep in mind, one that motivated much of the development of the Theory in the 1980s, is the \emph{Calder\'on Commutator} defined as follows.  
For Lipschitz function $ A $ on $ \mathbb R $, let 
\begin{equation*}
C _{A} f (x) := \int  f (y) \frac {A (x) - A (y)} {(x-y) ^2 } \; dy \,.  
\end{equation*}
Note that we have $ C_A = [M_A , \frac d {dx}] H $, where $ M_A$ is the operation of multiplication by $ A$.   We require $ A$ to be Lipschitz, as then we have, in some average sense, $\frac {A (x) - A (y)} {(x-y) ^2 } \simeq (x-y) ^{-1}  $.  And the deep fact is that we have 
$ \lVert C _{A}\rVert_{2 \mapsto 2} \lesssim \lVert A\rVert_{\textup{Lip}}$.  But, this is not at all easy to prove! .

A general definition of  Calder\'on-Zygmund operators we will consider.  

\begin{czo}\label{d.CZ} Let $0< \delta <1$ and 
Let $ K (x,y) \;:\; \mathbb R ^{d} \times \mathbb R ^{d} \backslash \{ (x,x) \;:\; x\in \mathbb R ^{d}\} 
\longrightarrow \mathbb R $ satisfy   kernel estimates 
\begin{gather} \label{e.s+s}
\lvert K (x,y)\rvert \le  C_T \lvert  x-y\rvert ^{-d-j}\,, \qquad x\neq y \in \mathbb R ^{d}\,, 
\\ \label{e.+s}
\lvert    K (x,y) -   K (x',y) \rvert  + 
\lvert    K (y,x) -   K (y,x') \rvert 
\le C_T \frac {\lvert  x-x'\rvert ^{\delta }}  {\lvert x-y\rvert ^{ d + \delta }} \,, 
\end{gather}
with the second condition holding provided $ \lvert x-x'\rvert < \tfrac 12 \lvert  x-y\rvert $.  Here, $0< C_T< \infty $. 
Occasionally, $ K (x,y)$ will be referred to as a Calder\'on-Zygmund kernel. 

Consider a linear operator  $ T: L^2\to L^2 $ such that 
\begin{equation*}
 T f (x) = \int  K (x,y) f (y) \; dy,\qquad x\notin\operatorname{supp}f,  
\end{equation*} 
for a fixed kernel $ K (x,y)$.

We then say that $ T$ is Calder\'on-Zygmund Operator, and write
$ T \in \operatorname {CZO} _{\delta } $ and  
\begin{equation}\label{e.czoNorm}
\lVert  T \rVert_{\textup{CZO} _ \delta  } \coloneqq 
\lVert T \rVert_{L ^2 (dx) \mapsto L ^{2} (dx)}  + C_ T   < \infty  \,. 
\end{equation}
\end{czo}

One should note that in one dimension, that the kernel $ K (x,y)= \frac 1 {\lvert  x-y\rvert }$ is a Calder\'on-Zygmund kernel, though the corresponding operator is not bounded. 
As well, it is hardly clear that the  Calder\'on Commutator is a bounded operator.  
Thus, it is a natural question to find a simple characterization of the  Calder\'on-Zygmund Operators. 
This class of operators was characterized by David and Journ\'e  \cite{MR763911}, in 
the famous $ T 1$ Theorem.

\begin{theorem}[$ T1$ Theorem]  An operator $ T$ with  Calder\'on-Zygmund kernel,  is $ L ^2 $-bounded if and only if  
for $ \mathbf T >0$,  we have the two uniform estimates  over all cubes $ I\subset \mathbb R ^{d}$.  
\begin{align}\label{e.T1}
\int _{I} \lvert T \chi _I\rvert ^2 dx \le \mathbf T ^2 \lvert  I\rvert \,, 
\\  \label{e.T*1}
\int _{I} \lvert T ^{\ast}  \chi _I\rvert ^2 dx \le \mathbf T ^2 \lvert  I\rvert \,. 
\end{align}
In the second line, $ T ^{\ast} $ is the adjoint of $ T$, namely it has the kernel $ K (y,x)$.  
\end{theorem}

The import of this result is that the full $ L ^2 $-inequality already follows from the boundedness of the operator on a very small set of functions, namely the indicators of cubes.  We should note that this is not the formulation of the Theorem as in \cite{MR763911},   but the version found in \cite{MR1232192}*{Chapter V}. 
Clearly, we prefer the form above, over its more familiar formulation, as it does not require the supplemental space $ BMO$. 
We refer to the two conditions \eqref{e.T1} and \eqref{e.T*1} as \emph{Sawyer testing conditions}, as their use in characterizing the bounded of operators first appeared in his two-weight Theorems on the Maximal Function \cite{MR676801} and the Fractional Integrals \cite{MR930072}. 

Let us use this Theorem to see that the Calder\'on Commutator is bounded.  Let us take the interval $ I=[a,b]$, and then using integration by parts, 
\begin{align*}
C_A (\chi _{I})(x) 
& = \int _{a} ^{b}  \frac {A (x) - A (y)} {(x-y) ^2 } \; dy   
\\
&= \frac {A (x) - A (b)} {(x-b) } - \frac {A (x) - A (a)} {(x-a) }+ 
\int _{a} ^{b} \frac {A' (y)} {x-y} \; dy \,. 
\end{align*}
The first two terms are bounded by $ \lVert A\rVert_{\textup{Lip}}$, and the third is the Hilbert transform applied to $ A' \chi _{(a,b)} \in L ^{\infty }$.  Hence the testing condition for $ C_A$ follows from the $ L ^2 $-boundedness of the Hilbert transform.  

\section{The $ A_p$ Weights } 

The $ A_p$ weights have the definition 

\begin{definition}\label{d.Ap}  For $ w $  an a.\thinspace e.\thinspace
positive function (a weight) on $ \mathbb R ^{d}$, 
 we define 
the $ A_p$ characteristic of $ w $ to be 
\begin{equation} \label{e.Ap}
\lVert w  \rVert_{ A_p}:= \sup _{I}  \lvert  I\rvert ^{-1}  \int _{I} w  \; dx  \cdot \Bigl[ 
\lvert  I\rvert ^{-1} \int _{I} w  ^{-1/(p-1)} \; dx\Bigr] ^{p-1} \,, \qquad 1<p<\infty \,, 
\end{equation}
where the supremum is over all cubes in $ \mathbb R ^{d}$.  In the case of $ p=1$, we set 
\begin{equation}\label{e.p=1}
\lVert w \rVert_{A_1} := \bigl\lVert \frac {M w} {w}\bigr\rVert_{\infty } \,. 
\end{equation}
\end{definition}

We note that $ \lVert w\rVert_{A_p}$ is not a norm, but continue to use the familiar notation.  Below, we will also write $ w (I) 
= \int _{I} w \; dx$ for the (non-negative) measure with density $ w$.  
It is a critical property, one that is key to the many beautiful properties of the $ A_p$ theory, that we have $ w >0$ a.e. 
In particular, this means that $ \sigma := w  ^{-1/(p-1)}$ is unambiguously defined.  Also, note that we have $ w \sigma ^{p-1} \equiv 1$, 
which casts the definition of $ A_p$ is a clear light: It requires that this pointwise equality continue to hold in an average sense, uniformly over all locations and scales.  

We will refer to $ \sigma $ as the \emph{dual measure.}
This language is justified by a useful observation from \cite{MR676801}.  The inequalities below are all equivalent for a 
linear operator $ T$: 
\begin{align}
\lVert T f\rVert_{ L ^{p} (w)}& \le C \lVert  f\rVert_{ L ^{p} (w)} \,, 
\\ \label{e.nat}
\lVert T (\sigma f)\rVert_{ L ^{p} (w)}&\le C \lVert T f\rVert_{ L ^{p} (\sigma )} \,, 
\\
\lVert T ^{\ast}  (w \phi )\rVert_{ L ^{p'} (\sigma )}&\le C \lVert T f\rVert_{ L ^{p'} (w) } \,.  
\end{align}
To pass from the first line to the second, use the change of variables $ f \mapsto \sigma \cdot f$.
There is a routine calculation, which is based on the basic identity of the weighted theory that $ p (p'-1)=1 $.  
And, note that the last line is the formal dual inequality to the second.  
Thus, the inequality \eqref{e.nat} expresses  duality in a natural way: Interchange the roles of $ w $ and $ \sigma $, and replace $ p$ by
dual index $ p'$.

Of course we are primarily interested in the case of $ p=2$.  
Two examples of $ A_2 $ weights to keep in mind, in dimension 1, are 
as follows.  First, for an arbitrary measurable set $ E\subset \mathbb R $, and $ N >0$,  the weight is $ w = N \chi _{E} + 
 \chi _{\mathbb R - E}$.   As long as $ \lvert  E\rvert>0 $, one has $ \lVert w\rVert_{A_2} \le \max(N,1/N)$.
 Indeed, we can assume $ N>1$. 
 For an interval $ I$ we have 
\begin{align}
\frac {w (I)} {\lvert  I\rvert}
\cdot  
\frac {\sigma  (I)} {\lvert  I\rvert}
& \le \frac {N \lvert  E\cap I\rvert + \lvert  E ^{c} \cap I\rvert  } {\lvert  I\rvert } \le N \,. 
\end{align}
This shows  that an $ A_2$ weight need not have any smoothness associated with it. 

A second example is the borderline case of $ w (x)= \lvert  x\rvert $.  This is \emph{not} an $ A_2$ weight as the dual measure $ \sigma (x) = \lvert  x\rvert ^{-1}  $ is not locally integrable.  But if we mollify the zero, setting for $ 0< \alpha  <1$, $ w _{\alpha } (x) = \lvert  x\rvert ^{\alpha } $, then we have $ \lVert w _{\alpha }\rVert_{ A_2} \simeq (1- \alpha  ) ^{-1} $.    
It is for such examples that one can verify that 
$
\lVert H \rVert_{ L ^2 (w _{\alpha })} \simeq (1- \alpha ) ^{-1} 
$. (Test on $ \chi _{[0,1]}$.)
But, these examples are somewhat misleading, in that the simple behavior of their zeros is not at all indicative of intricacy of the general $ A_2$ measure.

We comment on classical Theorem of Muckenhoupt \cite{MR0293384} concerning the $ A_p$ weights and the Maximal Function, defined by 
\begin{equation*}
M f (x) := \sup _{t>0} (2t) ^{-1} \int _{[-t,t] ^{d}} \lvert  f (x-u)\rvert  \; dt \,. 
\end{equation*}

\begin{theorem}\label{t.M} For $ w>0$, we have the following equivalences: 
\begin{enumerate}
\item  $ w\in A_p$; 
\item  $ M$ is bounded as a map from $ L ^{p} (w)$ to $ L ^{p,\infty } (w)$; 
\item  $ M$ is bounded as a map from $ L ^{p} (w)$ to $ L ^{p,\infty } (w)$. 
\end{enumerate}
Note that the weak and strong type norms are equivalent. 
\end{theorem}

Clearly, the strong-type inequality implies the weak-type. Using the formulation \eqref{e.nat}, and applying the maximal function to a the indicator of a cube directly proves that $ w\in A_p$.  
So, the content of the result is that the $ A_p$ property implies the strong type inequality.  
Here, the fact that $ w >0 $ a.e. is decisive, and the shortest--six lines--proof of this is due to Lerner \cite{MR2399047}.  
Nevertheless, it seems confusing that the weak and strong types should be equivalent.  
The sharp dependence of the Maximal Function on the $ A_p$ characteristic is helpful here.  
For $ w \in A_p$, how does the norm depend upon $ \lVert w\rVert_{A_p}$?  
Buckley \cite{MR1124164} studied the question and proved  

\begin{theorem}\label{t.buckley} For $ w\in A_p$, we have 
\begin{align}
\lVert M \rVert_{L ^{p} (w) \mapsto L ^{p,\infty } (w)} & \lesssim \lVert w\rVert_{A_p} ^{1/p} \,,
\\ \label{e.mapStrong}
\lVert M \rVert_{L ^{p} (w) \mapsto L ^{p } (w)} & \lesssim \lVert w\rVert_{A_p} ^{1/(p-1)}  
\end{align}
\end{theorem}

Thus, the norm dependence is rather different.  This is a basic set of inequalities, due to the notion of Rubio de Francia extrapolation 
\cite{MR745140}.

A final, critical property for us is the so-called $ A _{\infty }$-property. It states that an $ A_p$ weight cannot be too concentrated in any one cube.  Indeed, as we will illustrate in the context of the Maximal Function, this is the single property of $ A_p$ weights that can be used to prove sharp results, and it can only be used once. 

\begin{lemma}\label{l.Ainfty}  Let $ w \in A_p$, $ I$ is a cube and $ E \subset I$.  We then have 
\begin{equation}\label{e.Ainfty}
\frac {\lvert  E\rvert } {\lvert  I\rvert }\le  \lVert w\rVert_{A_p} ^{1/p} 
\Bigl[\frac {w (E)} {w (I)} \Bigr] ^{1/p} \,. 
\end{equation}
\end{lemma}
    
\begin{proof} 
The property that $ w >0$ a.e. allows us to write 
\begin{align*}
\frac {\lvert E\rvert} {\lvert I\rvert}  & = 
\frac {\int _{E} w ^{1/p} (x) w (x) ^{-1/p} \; dx } {\lvert  I\rvert } 
\\
& \le  \frac {w (E) ^{1/p}  \sigma (I) ^{1/p'}} {\lvert  I\rvert } 
\\
& = \Bigl[  \frac {w (E) } {w (I)} \Bigr] ^{1/p}   \frac {w (I) ^{1/p} \sigma (I) ^{1/p'}} {\lvert  I\rvert } 
\end{align*}
which proves the Lemma. 
\end{proof}

\section{Dyadic Grids} 

Combinatorial arguments, stopping time arguments or decompositions of functions and operators, will frequently be done with the help of \emph{dyadic grids}.  In this section, we collect a number of elementary facts that we will need from time to time.  
At different moments, the methods and constructions of this section will in fac be decisive for us.  

By a \emph{grid} we mean a collection $ \mathcal G$ of cubes in $ \mathbb R ^{d}$ with  $ I \cap I' \in \{\emptyset , I, I'\}$ for all $ I, I' \in \mathcal I$.   
The cubes can be taken to be a product of clopen intervals, although the behavior of functions or weights on on the boundary of cubes in a grid will never be a concern for us.  If $ G, G' \in \mathcal G$, with $ G'$ the smallest element of $ \mathcal G$ that strictly contains $ G$, we refer to $ G'$ as the \emph{$ \mathcal G$-parent} of $ G$, and $ G'$ is a $ \mathcal G$-child of $ G'$.  Let $ \textup{Child} _{\mathcal G} (G')$ denote the collection of all  $ \mathcal G$-children of $ G'$.  If the grid is understood, the $ \mathcal G$ is suppressed in the notation. 

We will say that $ \mathcal G$ is a  \emph{dyadic grid} if each cube $ I \in \mathcal G$ these two properties hold. 
(1)  $ I$  is the union of $ 2 ^{d}$-subcubes of  equal volume (the children of $ I$), 
and (2) the set of cubes $ \{ I' \in \mathcal G \;:\; \lvert  I'\rvert= \lvert  I\rvert  \}$ partition $ \mathbb R ^{d}$.

Associated to any dyadic grid $ \mathcal D$ are the usual conditional expectations and martingale differences are given by 
\begin{align}\label{e.expect}
\mathbb E _{I} f & := \chi _{I} \cdot  \lvert  I\rvert ^{-1} \int _{I} f \; dx  \,, 
\qquad 
\Delta _{I} f  := \sum_{I' \in \textup{Child} (I)} \mathbb E _{I'} f  - \mathbb E _{I} f \,. 
\end{align}
And, we also set 
\begin{align*}
\mathbb E _{k} f
:= \sum_{\substack{I\in \mathcal D\\ \ell (I)= 2 ^{-k} }} \mathbb E _{I} f \,, 
\qquad 
\Delta _{k} f := \sum_{\substack{I\in \mathcal D\\ \ell (I)= 2 ^{-k} }} \Delta_{I} f  \,. 
\end{align*}
Then, by the Martingale Convergence Theorem, for $ f\in L ^{1} (dx)$, $ \mathbb E _{k} f \rightarrow f $ a.e. And, by the Muckenhoupt Theorem for the Maximal Function,  for $ w\in A_2$, and $ f\in L ^2 (w)$, the same conclusion holds.

\subsection{Proof of Buckley's Maximal Function Inequality}

As an illustration of the use the $ A _{\infty }$ condition, let us return to Buckley's estimate \eqref{e.mapStrong}, 
and prove it in the dyadic case.  Namely, for choice of dyadic grid $ \mathcal D$ in $ \mathbb R ^{d}$, we define the 
associated Maximal Function 
\begin{equation} \label{e.Mdy}
M f (x) := \sup _{I\in \mathcal D}  \chi _{I} (x) \mathbb E _{I} \lvert  f\rvert  
\end{equation}
where here we have introduced the notation $ \mathbb E _{I} \phi := \lvert  I\rvert ^{-1} \int _{I} \phi  $.  
Also, we are continuing with the same notation for the Maximal Function, suppressing its dependence on the choice of grid. 
For this operator, we will prove \eqref{e.mapStrong}.

We make the definition of the stopping cubes.  

\begin{definition}\label{d.1stop} Let $ \mathcal G$ be a grid, $ \sigma $ a weight.  Given cube $ I \in \mathcal G$, 
we set the \emph{stopping children of $ I_0$, written $ \mathcal C (I)$,} to be the maximal dyadic cubes $ I'\subset I$
for which $  \mathbb E _{I'} \sigma > 4 \mathbb E _{I} \sigma $.    A basic property of this collection is that 
\begin{equation}\label{e.14}
\sum_{I'\in \mathcal C (I)} \lvert  I'\rvert< \tfrac 1 4 \lvert  I\rvert\,.   
\end{equation}

We set the \emph{stopping cubes of $ I$} 
to be the collection $ \mathcal S (I) = \bigcup _{j\ge 0} \mathcal S_j (I)$, where we inductively define $ S_0 (I) := \{I\}$, 
and $ S _{j+1} (I) = \bigcup _{I'\in \mathcal S _{j} (I)} \mathcal C (I)$.  Thus, these are the maximal dyadic cubes, so that 
passing from parent to child in $ \mathcal S$, the average value of $ \sigma $ is increasing by at least factor $ 4$.  
\end{definition}

\begin{proof}[Proof of \eqref{e.mapStrong}.]   It is the fundamental Theorem of Eric Sawyer \cite{MR676801} that for the Maximal Function, we have a powerful variant of the David Journ\'e $ T1$ Theorem.  Namely, for \emph{any} pairs of weights $ (w, \sigma )$, we have the equivalence between these two inequalities
\begin{align}
\lVert M (\sigma f)\rVert_{L ^{p}(w)} &\le C_1 \lVert f\rVert_{L ^{p}(\sigma) }  \,,
\\ \label{e.sawT}
\int _{I} M (\sigma \chi _{I}) ^{p} w \; dx &\le C_2 ^{p} \sigma (I)\,, \qquad I\in \mathcal D \,. 
\end{align}
Moreover, letting $ C_1$ and $ C_2$ be the optimal constants in these two inequalities, we have $ C_1 \simeq C_2$.   
Notice that this shows that the Maximal Function bound reduces to a testing condition.  

And so, in the special case that $ w \in A_p$, and $ \sigma = w ^{1-p'}$, we estimate the constant $ C_2$.  
\begin{equation}\label{e.zm2}
\int _{I} M   ( \sigma \chi _ I) ^{p} \; w (dx) \lesssim \lVert w\rVert_{A_p} ^{p'} \sigma (I) \,. 
\end{equation}
We do so by passing to the stopping cubes $ \mathcal S (I)$, and estimating as below, where we will use some common manipulations in the $ A_p $ theory. 
\begin{align}
\int _{I} M  ( \sigma \chi _I) ^{p} \; w (dx) 
& \le 
\int _{I} \Bigl[ \sum_{S \in \mathcal S (I)}  \frac {\sigma (S)} {\lvert  S\rvert } \cdot \chi _ S \Bigr]
 ^{p} \; w (dx) 
 \\  \label{e.gtrk}
 & \lesssim 
 \sum_{S \in \mathcal S (I)} \int _{I}\Bigl[ \frac {\sigma (S)} {\lvert S\rvert } \cdot \chi _S \Bigr]
 ^{p} w (I) 
 \\  \label{e.swp1}
 & \le \lVert w\rVert_{A_p}   \sum_{S \in \mathcal S (I)} \sigma (S) 
\\   \label{e.ainfy1}
& \lesssim    \lVert w\rVert_{A_p}   \lVert \sigma \rVert_{A_{p'}} \sigma (I) 
\\  \label{e.apid1}
& =  \lVert w\rVert_{A_p}  ^{p'} \sigma (I) \,. 
 \end{align}
 This proves our estimate.  
Here, we have taken these steps.  
\begin{description}
\item[\eqref{e.gtrk}]  Pointwise, the sum $  \sum_{S \in \mathcal S (I)}  \frac {\sigma (S)} {\lvert  S\rvert } \cdot  \chi _S (x)$ is 
super-geometric, so comparable to its maximal term in the summand. This allows us to move the $ p$th power inside the sum. 
\item[\eqref{e.swp1}]  We are using the definition of $ A_p$ here.  
\item[\eqref{e.ainfy1}] The $ A_ \infty $ property is decisive. By \eqref{e.14} and \eqref{e.Ainfty}, we have, using the notation for the 
stopping children from Definition~\ref{d.1stop}, 
$ \sum_{I' \in \mathcal C (I) } \sigma (I) \le (1- c \lVert \sigma \rVert_{A _{p'}} ^{-1} ) \sigma (I)$, permitting us to sum a geometric series 
to get this estimate. 
\item[\eqref{e.apid1}]  By inspection, $  \lVert \sigma \rVert_{A_{p'}} =  \lVert w\rVert_{A_p}  ^{p'-1}$.  
\end{description}
\end{proof}

\subsection{Random Dyadic Grids, Good and Bad Cubes}\label{s.good}

We are used to thinking of a dyadic grid as being canonical, namely the cubes  
\begin{equation*} \mathcal D := 
\bigl\{  2 ^{k} ( n + [0,1) ^{d}) \;:\; k\in \mathbb Z\,,\ n \in \mathbb Z ^{d}  \bigr\} \,. 
\end{equation*}
This choice has a strong \emph{edge effect}, it for instance distinguishes the origin, in that it is the vertex  of infinitely many cubes.  This sort of anomaly on the other hand should be typically rare. 
Quantifying this is achieved by a random grid.  To present one typical example, if $\mathcal G$ is a dyadic grid, and the interval $ [0,1) ^{d}$ is in $ \mathcal G$, it has one $ 2 ^{d}$ possible parents, found by taking the cube to be the product of one of the 
two intervals $ [0,2)$ and $ [-1,1)$ in each coordinate separately.   To randomize $ \mathcal G$, these possible choices of grids should be equally likely.

 For any $\beta =\{\beta
_{l}\}\in   \boldsymbol \beta :=  \bigl\{\{0,1\}^{d} \bigr\}^{\mathbb Z }$,  and cube $ I$, set 
\begin{equation} \label{e.iI}
I \dot + \beta = I + \sum_{ l < \ell (I)} \beta _{l} 2 ^{-l}\,, 
\end{equation}
where $ \ell (I):= \lvert  I\rvert ^{1/d}  $  is the side length of the cube.
Then, define the dyadic grid ${\mathbb{D}}_{\beta }$ to be the collection of cubes 
$
{\mathbb{D}}_{\beta }= \{I\dot + \beta \;:\; I \in \mathcal D\}
$.  
This parametrization of dyadic grids appears explicitly in \cite{MR2464252},
and implicitly in \cite{MR1998349}*{section 9.1}. 

Place the uniform  probability measure $\mathbb{P}$ on the space $ \boldsymbol \beta $.  
Namely, the probability that any coordinate $ \beta _j $ takes any one value in $ \{0,1\} ^{d}$ is $ 2 ^{-d}$, and the coordinates $ \beta _j$
are independent of one another.  

Let us see how the randomization affects the edge effect mentioned above.  Let $ 0 < \gamma <1$ be a fixed parameter, and $ r \in \mathbb Z _+$ is a fixed integer.   We say that   say that a pair of intervals $ (I, J) \in \mathcal D _{\beta }$ are \emph{good} if 
the smaller interval, say  $ I$, satisfies 
$ 2 ^{r}\ell (I) < \ell (J)$, and  
\begin{equation}\label{e.good}
\textup{dist} (I,\partial J) \ge  \ell (I) ^{\gamma } \ell (J) ^{1- \gamma }\,. 
 \end{equation}
And an interval $ I$ is said to be \emph{good}, if for all intervals $ J$ with $ \ell (J)> 2 ^{r} \ell (I)$, we have that the pair $ (I,J) $ is good.  
Otherwise, we say that the cube is \emph{bad}. 

An important property of goodness is the independence of the location or scale of a cube $ I$ and its goodness.  Take $ I \dot+ \beta \in \mathcal D _{\beta }$. 
The spatial position of $ I$ is given by the formula \eqref{e.iI}, 
which only depends upon $ \beta _j$ for $  2 ^{-j} < \ell (I)$.  And, for a larger cube $ J$, the position of $ J$ can be written as 
\begin{equation*}
J+ \sum_{j \;:\; 2 ^{-j} < \ell (I)} 2 ^{-j} \beta _j  + \sum_{j \;:\; \ell (I) \le 2 ^{-j} < \ell (J)}  2 ^{-j} \beta _j \,. 
\end{equation*}
And hence, the position of $ J$ \emph{relative to $ I$} depends only on the coordinates $ \beta _j$ for $  \ell (I) \le 2 ^{-j} < \ell (J)$, 
and hence is independent of the location of $ I$.

As a consequence the probability of a given cube is bad is independent of the location or scale of $ I$. Denoting this  probability by 
$ \pi _{r, \gamma }$, it is an elementary exercise to see that $ \pi _{r, \gamma } \lesssim 2 ^{-r \gamma }$.  As it will turn out, it will be sufficient to have this probability less than one, for a choice of $ \gamma $ that depends upon the Calder\'on-Zygmund Operator $ T$, 
and can be taken to be a small multiple of the constant $\delta  $ in the Definition~\ref{d.CZ}. 

\subsection{Haar Shifts, Dyadic Calder\'on-Zygmund Operators}

In one dimension, the Martingale Difference in \eqref{e.expect} is given by the rank-one projection 
$
\Delta _{I} f =  \langle f, h _{I} \rangle \cdot h_I
$
where $ h_I$ is the  Haar function, given by $ h _I := (- \chi _{I_-} + \chi _{I_+}) \lvert  I\rvert ^{-1/2} $, where $ I _{\pm}$ denotes the two children of $ I$.  And then, the simplest possible dyadic Calder\'on-Zygmund operator would be a \emph{martingale transform}
\begin{equation*}
T f := \sum_{I} \varepsilon _I   \langle f, h _{I} \rangle \cdot h_I \,. 
\end{equation*}
The amenability of these operators to issues of measurability, and stopping time arguments has long been exploited, leading to a 
remarkable set of properties that are known for these objects.

Below, we will say that martingale transforms have \emph{complexity $1$}.  To motivate this upcoming definition, let us recall the remarkable result of Stephanie Petermichl, concerning the Hilbert transform.  In one dimension, consider the \emph{dual} to the classical Haar function given by $ g_I = (-h _{I_-} + h _{I_+})/\sqrt 2$, and the special operator given by 
\begin{equation*}
U f = U _{\beta } f := \sum _{I\in \mathcal D _{\beta }} \langle f, g_I \rangle \cdot h_I \,. 
\end{equation*}
The Hilbert transform can be recovered from the operators $ U _{\beta }$, namely the result below holds. 

\begin{theorem}\label{t.sp}  Let $ \textup{Dil} _{\delta } f (x) = f (x/\delta )$.  For non-zero constant $ c$, we have 
\begin{equation*}
\mathbb E _{\beta} \int _{1} ^{2}   \textup{Dil} _{\delta } U _{\beta }  \textup{Dil} _{1/\delta } \; \frac {d \delta } \delta = c H  
\end{equation*}
Here, the expectation is taken over $ \beta \in \boldsymbol \beta $.  
\end{theorem}

The Hilbert transform is distinguished by different properties, including being $ L ^2 $-bounded, translation and dilation invariant, and (formally) satisfying $ H (\cos)=c \cdot \sin $.  By inspection, $ U _{\beta }$ is $ L ^2 $-bounded.  The averaging procedure above provides translation invariance, and dilation invariance, as we have used the Haar measure for the dilation group in the average.  For the last property, note that $ g_I$ is a localized cosine, while $ h$ is a localized sinus.  We refer the reader to \cites{MR1756958,MR2464252} for a precise proof of this Theorem. 

The import of this result is that in situations where there is a translational and dilational invariance, one can prove results about the Hilbert transform by considering the much simpler operators $ U$---where \emph{tail behavior is no longer an issue}.  Similar representations are available for other distinguished convolution kernels. For instance, the Beurling operator \cite{MR1992955} can be recovered from martingale transforms, 
while the Riesz transforms are closer to the Hilbert transform \cite{MR1964822}.  The most general result known in this direction is 
\cite{0911.4968}, which shows that all smooth, odd one dimensional Calder\'on-Zygmund kernels can be obtained by a variant of Stephanie Petermichl's method.  

  A more general definition is as follows. 
In higher dimensions, we mention that the martingale differences are finite rank projections, but there is no canonical choice of the Haar functions in this case.  Below, by \emph{Haar function} we will a function $ h _{I}$, supported on $ I$, constant on its children, and orthogonal to $ \chi _{I}$ (and no assumption on normalizations).   
And, by a  \emph{generalized Haar function}  as a function $ h_I$ which is a linear combination of $ \chi _I$, and $ \{\chi _{I' } \;:\; I'\in \textup{Child} (I)\}$. 
Such a function supported on $ I$ but need not be orthogonal to constants.  

\begin{definition}\label{d.} For integers $(m,n) \in \mathbb Z _+ ^2 $, we say that  linear operator $ S$ is a \emph{ (generalized) Haar shift operator of parameters $ (m,n)$} if 
\begin{equation}\label{e.mn}
S f (x) = \sum_{I \in \mathcal D}\; \sideset {} { ^ {(m,n)}} \sum_{\substack{I',J'\in \mathcal D\\ I',J'\subset I }} 
\frac { \langle f, h ^{I'} _{J'} \rangle} {\lvert  I\rvert } h ^{J'} _{I'} 
\end{equation}
where (1) in the second sum, the superscript $ ^{(m,n)}$ on the sum means that in addition we require $ \ell (I') = 2 ^{-m} \ell (I)$ and $ \ell (J')= 2 ^{-n} \ell (I)$, and (2) the function $  h ^{I'} _{J'}$ is a (generalized) Haar function on $ J'$, and $  h ^{J'} _{I'}$ is one on $ I'$, with the joint normalization that 
\begin{equation} \label{e.normal}
\lVert  h ^{I'} _{J'}\rVert_{\infty } \lVert  h ^{J'} _{I'}\rVert_{\infty } \le 1 \,. 
\end{equation}
In particular, this means that we have the representation 
\begin{equation} \label{e.NOR}
S f (x) = \sum_{I\in \mathcal D} \lvert  I\rvert ^{-1} \int _{I} f (y)  s _{I} (x,y) \; dy  
\end{equation}
where $ s_I (x,y)$ is supported on $ I \times I$, with $ L ^{\infty }$ norm at most one.  
We say that the \emph{complexity} of $ S$ is $ \max (m,n)$.  
\end{definition}

These are dyadic variants of Calder\'on-Zygmund operators.  Note in particular that \eqref{e.normal} is analogous to \eqref{e.s+s}, while the `smoothness' criteria is replaced by the parameters $(m,n)$.   
Consider a Haar shift operator. It is an $ L ^2 $-bounded operator, in particular its norm is at most one.  The situation for generalized shifts is far more subtle, and here, we should single out the following definition, for it distinguished role in the theory, though not necessarily this paper. 
  We call an operator $ S$ a paraproduct if it is a generalized Haar shift of parameters $ (0,1)$ or $ (1,0)$.  
To be specific, it, or its dual, is of the form 
\begin{equation} \label{e.para}
S f = \sum_{I\in \mathcal D} \mathbb E _{I} f  \cdot h _{I}
\end{equation}
where $ h_I$ is a Haar function.   A fundamental fact here is the following special case of the $ T1$ Theorem, in the dyadic case. 

\begin{theorem} Let $ S$ be as in \eqref{e.para}.  Then, $ S$ is $ L ^{2}$-bounded if and only if we have 
\begin{equation} \label{e.paraTesting}
\lVert S \chi _{I}\rVert_{2} \lesssim \lvert I\rvert ^{1/2} \,. 
\end{equation}
\end{theorem}

This is a particular variant of the famous Carleson Embedding Theorem, and the main step in extending the David Journ\'e $ T1$ Theorem to the dyadic setting.  

More generally, we have the following quantitative form of the Dyadic $ T1$ Theorem.  

\begin{theorem}\label{t.dT1} Let $ S$ be a generalized Haar shift operator of complexity $ \mu $.  Then $ S$ extends to a bounded operator 
on $ L ^{2} (\mathbb R ^{d})$ if and only if we have, uniformly over cubes $ I$, 
\begin{align}
\int _{I} \lvert  S \chi _{I}\rvert ^2 \; dx &\le \mathbf S ^2 \lvert  I\rvert \,,   
\\
\int _{I} \lvert  S ^{\ast}  \chi _{I}\rvert ^2 \; dx &\le \mathbf S ^2 \lvert  I\rvert \,,   
\end{align}
Moreover, we have $ \lVert S\rVert_{ 2 \mapsto 2} \lesssim \mu \mathbf S + \mu ^2 $. 
\end{theorem}

There are two points to make here. The first is that there is a  weak dependence of the norm of the operator as a function of the complexity $ \mu $.  The second, is the familiar, but not mentioned to this point, feature of the Calder\'on-Zygmund theory, that thee operators have  strong features.   If $ S$ is a bounded operator, then, the sum in \eqref{e.NOR} is unconditional in $ I$.  
The import of this feature, important for proof of the main Theorem, is that  decompositions of dyadic cubes lead immediately to  decompositions of operators. 
In the second, an $ L ^2 $-bounded Calder\'on-Zygmund operator is necessarily bounded on many other spaces. 
Of particular interest for us is the endpoint estimate for $  L ^{1}$: 

\begin{theorem}\label{t.L1} Let $ S$ be a dyadic shift operator of complexity $ \mu $,  which is bounded on $ L ^2 (\mathbb R ^{d})$. Then, we have the estimate 
\begin{equation}\label{e.L1}
\sup _{\lambda >0} \lambda \lvert  \{ S f > \lambda \}\rvert \lesssim 
\{(1 + \lVert S\rVert_{ 2 \mapsto 2}) ^2 + \mu \} \lVert f\rVert_{1} \,. 
\end{equation}
\end{theorem}

This is a well-known principle, but the weak-dependence on the complexity is a point observed by Hyt\"onen. See \cite{1010.0755}*{Theorem 5.2}.

\subsection{A Weighted Version of the Dyadic $ T1$ Theorem}

A crucial step is to prove a weighted version of the $ T1$ Theorem, one that holds for general weights.  To emphasize this point, for a pair of weights $ (w, \sigma )$, which are not necessarily related, we set the \emph{two weight $ A_2$ } condition to be 
\begin{equation}\label{e.2a2}
\lVert w, \sigma \rVert_{A_2} := \sup _{I\in \mathcal D} \frac { w (I)} {\lvert  I\rvert }  \frac { \sigma (I)} {\lvert  I\rvert }  \,. 
\end{equation}
We have this variant of the $ T1$ Theorem, for generalized Haar shift operators, in the weighted setting. 

\begin{theorem}\label{t.well}  Let $ S$ be a generalized Haar shift operator of complexity $ \mu $, and $ (w ,\sigma )$ a pair of weights. 
We have $ \lVert S (\sigma f)\rVert_{L ^2 (w)} \le C \lVert f\rVert_{\sigma }$, where 
\begin{gather} \label{e.C<}
C \lesssim _{d} \mu \mathbf S  + \mu ^2 \lVert w, \sigma \rVert_{A_2} ^{1/2} 
\\
\int _{I}  S (\sigma \chi _I) ^2  \; w (dx)  \le \mathbf S ^2 \sigma (I) \,, 
\\
\int _{I}  S ^{\ast}  (w  \chi _I) ^2  \; \sigma  (dx)  \le \mathbf S ^2 w  (I) \,, 
\end{gather}
\end{theorem}

Here, we are considering the weighted inequality in its natural form, see \eqref{e.nat}.  
And we are bounding the weighted norm of the Haar shift in terms of the two-weight $ A_2$ condition, as well as the testing condition. 
Of particular importance for the proof of the linear bound is the \emph{very weak dependence} of the constants on the $ A_2$ characteristic.  For the proof, see \cite{MR2407233} and for the quantitative estimate above \cite{1010.0755}*{Theorem 3.4}.  
In particular, the proof is a weighted variant of the usual proof of the dyadic $ T1$ Theorem, with an important point being that one should use \emph{weighted Haar functions} to give the proof.

\section{The Random BCR Algorithm} 

A proof of the $ T1$ Theorem must, implicitly, or explicitly, decompose the Calder\'on-Zygmund operator into appropriate components. 
In the language of random dyadic shifts, the remarkable result of \cite{1010.0755}*{Theorem 4.1} is 

\begin{theorem}\label{t.bcr} Let $ T$ be a Calder\'on-Zygmund Operator $ T$ with smoothness parameter $ \delta $.  Then, we can write 
\begin{equation} \label{e.bcr}
T = C  \mathbb E _{\beta } \sum_{ (m,n)\in \mathbb Z _+ ^2 }  2 ^{- (m+n) \delta /2} S _{m,n} ^{\beta } 
\end{equation}
where (a) the expectation is taken over the space of random dyadic grid; $ S _{m,n}$ is a (random) dyadic shift; (c) the shifts of parameters $ (0,1)$ and $ (1,0)$ are generalized shifts; (d) all other shifts need not be generalized; (e) the constant $ C$ is a function of $ T$, and the smoothness parameter $ \delta $.   In particular, we will have, uniformly over the probability space, 
\begin{equation*}
\lVert S _{m,n} ^{\beta }\rVert_{2 \mapsto 2} \lesssim 1 \,. 
\end{equation*}
\end{theorem}

The focus with Theorem~\ref{t.sp} is noteworthy. The prior result obtains the Hilbert transform as a convex combination of Haar shifts of \emph{bounded} complexity. The Theorem above obtains it as a  sum of Haar shifts, but one that is rapidly converging in complexity.

In the dyadic setting, similar results were proved by Figiel \cite{MR1110189}, and independently by \cite{MR1085827}, with the latter article being broadly influential. 
The method of expanding Calder\'on-Zygmund operators using this method reveals subtle approximation theory properties of these operators. 
This method is not random, but has the disadvantage of using operators which are not purely dyadic.  

Indeed, the Theorem above looks wrong. Using standard Haar basis in one dimension, the inner product $ \langle H h _{[0,1]}, h _{[0, 2 ^{k})} \rangle$ does not have the good decay properties in terms of complexity claimed above. Instead, one needs a concept like the goodness property of \S\ref{s.good}. 
And indeed, this is the main point, the inner product $ \langle H h _{I}, h _{J} \rangle$ will be small, if the pair of intervals $ (I,J)$ are good. 

\section{The Corona and the Linear Bound} 

In the prior proofs of the linear bound for operators, one used the averaging technique of Petermichl, see Theorem~\ref{t.sp}, to represent the Calder\'on-Zygmund operator as an average of Haar shifts of \emph{bounded} complexity.  And then, verified the linear bound for such shifts. 
But, the representation \eqref{e.bcr} gives one another option.  
For an $ A_2$ weight, and an arbitrary Haar shift operator $ S$, verify the linear bound, with only moderate growth in the complexity $ \mu $ of the Haar shift.  Here, we can allow any polynomial dependence on the complexity.  
We have already described this in two different places, the first is the 
dyadic two-weight $ T1$ Theorem, Theorem~\ref{t.well}, and the second is the weak-$L^1$ inequality, \eqref{e.L1}.  

The relevant result is \cite{1010.0755}*{Equation (5.5)}. 

\begin{theorem}\label{t.testing} Let $  S$ be a generalized Haar shift operator of complexity $ \mu $, and $ \mathbf S=\lVert S\rVert_{2 \mapsto 2} $.  For $ w \in A_2$ and $ \sigma = w ^{-1} $, and cube $ I$, we have 
\begin{equation*}
\Bigl[\int _{I} \lvert S \sigma \chi _I\rvert ^2 \; w (dx) \Bigr] ^{1/2}  \lesssim  ( \mu +1 ) (\mathbf S+ \mu +1 ) \lVert w \rVert_{A_2} 
\sigma (I) ^{1/2} \,. 
\end{equation*}
\end{theorem}

The method of proof here, aside from the dependence on the complexity, is derived from \cite{MR2657437}, and is a subtle extension of the method  used  to prove \eqref{e.mapStrong}, the sharp dependence on the $ A_p$ characteristic for the Maximal Function.  
Indeed, the interested reader should first consult \cite{MR2657437}, which does not seek to track dependence of the bound on the constants. 
 This argument  uses the stopping cubes, as given in Definition~\ref{d.1stop}. And, this decomposition is then used to decompose the operator.
Then, the main step is to  identifies a notable extension of the John-Nirenberg inequalities that holds in the two-weight setting, for the decomposed operator.  With this, we conclude our discussion of the proof of the linear bound for Calder\'on-Zygmund operators.

\section{History} 

\paragraph{}  The weighted theory came of age with the paper \cite{MR0312139} of Hunt--Muckenhoupt--Wheeden, showing that for non-negative weight $ w$, the Hilbert transform is bounded on $ L ^2  (w)$ if and only if $ w\in A_2$.  
Still, early proofs combined properties of the weight, including the $ A_ \infty $ property we have used, with the Reverse Holder inequality, and the good-lambda technique, to deliver estimates for the norm of the Hilbert transform of the order of $ \lVert w\rVert_{A_2} ^2 $.  
These and the other comments about history reflect the authors' knowledge, but as he was not a participant in the development of the subject, they will certainly be incomplete. Apologies for omissions and gaps are extended in advance. 

\paragraph {}  The rapid development of the $ A_p$ theory in the 1970's lent some credence to the thought that similar  variants of the $ A_p$ condition could be used to characterize the two-weight inequalities as well.   The characterization for the  Hardy operator \cite{MR0311856} confirmed this.  
It was a surprise when Sawyer \cite{MR676801} showed that such conditions cannot be used for the Maximal Function, instead one must use the testing conditions in \eqref{e.sawT}. (For a little more detail, consult the counterexample discussed in Sawyer's paper.)   

\paragraph {}  In the two-weight setting, the Hardy operator is somehow the easiest to study, the Maximal Function is the next step harder, then the fractional integrals, and finally the singular integrals.   It took several years for the proof of the two-weight inequalities for the fractional integrals to be characterized.  Sawyer gave the characterization in the $ T1$ language in  \cite{MR930072}. 
This was contemporaneous with the David-Journ\'e $ T1$ Theorem, but the connection was not widely appreciated until much later, especially by the work of Nazarov-Treil-Volberg. For history on this last point, see \cite{V}. 

\paragraph {}  In the two-weight setting, one can have the fractional integral operators mapping $ L ^{p}$ into $ L ^{p}$, indeed this is the hard case.  In the case of $ L ^{p}$ being mapped into $ L ^{q}$, for $ q>p$, there is a second characterization due to \cite{gk}, also see 
\cite{MR1791462}*{Chapter 3}, and \cite{MR1175693}.  This characterization can be used to prove the sharp $ A _{p,q}$ bound for the fractional integrals on $ \mathbb R ^{d}$, see \cite{MR2652182}.   

\paragraph {} The paper of Sawyer-Wheeden \cite{MR1175693} extends the two-weight inequality for the fractional integrals to homogeneous spaces; this is an interesting direction, which has been, and will be, explored in many different directions.  

\paragraph {} The question of the sharp dependence of the norm estimates of different operators, in terms of the $ A_p$ characteristic was specifically raised by Buckley \cite{MR1124164}, where the estimate \eqref{e.mapStrong} was proved.   These bounds for the Maximal Function, together with the Rubio de Francia extrapolation technique leads to an important simplification of the analysis of many of the weighted inequalities.  Namely, as is demonstrated in \cite{MR2140200}, identifying a sharp exponent in $ A_p$ characteristic for a single distinguished choice of $ p$ can prove the entire range of inequalities.  For the Calder\'on-Zygmund operators, this index is $ p=2$.  

\paragraph {} 
 In a different direction, Fefferman and Pipher \cite{MR1439553} recognized the interest of this question, for singular integrals, with the weights  $ w \in A_1$.  Wittwer  \cite{MR1748283} proved the linear bound for $ A_2$ weights, for martingale transforms.  
 Petermichl and Volberg \cite{MR1894362} showed the same for the Beurling operator, proving a conjecture of Astala on quasi-conformal maps as a consequence. 
 Much later, a certain two-weight inequality for the Beurling operator was proved \cite{0805.4711} as a crucial step in proving another conjecture  of Astala.  These examples motivate in part the interest in such questions.  
 Other motivations are derived from considerations in spectral theory \cite{kwon-2007}, operator theory \cite{NV}, and orthogonal polynomials \cite{MR2465289}.

 \paragraph {} It was an important breakthrough when Stefanie Petermichl proved the linear bound for the Hilbert transform \cite{MR2354322}. This technique was based on the one hand, the representation of the Hilbert transform as an operator of complexity one, and on the other   on the Bellman function method. 
The latter, deep, technique could require substantive modification if the Haar shift changes; these modifications were spelled out for the Riesz transforms in \cite{MR2367098}, and dyadic paraproducts \cite{MR2433959}.  

\paragraph {} An inequality used in some of these developments was the so-called bilinear embedding inequality of Nazarov-Treil-Volberg, \cite{MR1685781}. 
The latter is a deep extension of the (weighted) Carleson embedding inequality to a two-weight setting.  
This inequality can also be interpreted in the language of fractional integrals, and the Sawyer method can be used to prove it, and extend it to other $ L ^{p}$ settings \cite{0911.3437}, as well as vector-valued settings \cite{1007.3089}.

\paragraph {} Andrei Lerner \cite{localSharp} devised a remarkable inequality, giving pointwise control of a function in terms of a sum of local oscillations.  This inequality can be used to provide equally remarkable proofs of the sharp $ A_p$ inequalities for dyadic Calder\'on-Zygmund operators \cites{MR2628851,1001.4254},  even in certain vector-valued situations.  
As of yet, it is not understood how to use this method on continuous Calder\'on-Zygmund operators.

\paragraph {} Commutators of the form $ [T, M_b]$ are of interest, for instance, the Calder\'on Commutator can be written in this form. 
And the paper  of Chung-Pereyra-P{\'e}rez  \cite{1002.2396} gives a  complete discussion of this question in the setting of $ A_p$ weights.  The two-weight variants appear to be largely open.

\paragraph {} Lerner conjectured that the Littlewood-Paley Square function would have a different behavior in terms of its $ A_p$ characteristic. Namely, the case of $ p=3$ was the critical index, and the power on the $ A_p$-characteristic was $ 1/2$.  
He used his `local oscillation' inequality, as well as other considerations, to prove this inequality in full generality \cite{1005.1422}. 

\paragraph {} The paper \cite{MR2657437} proved the $ A_2$ linear bound for all Haar shifts, using a Corona decomposition that has been useful to the complete resolution of the Conjecture.  
The technique is obtaining a natural Corona decomposition in order to verify the testing conditions. 
This paper gave a rather poor dependence in terms of the complexity of the Haar shift parameter, but the role of complexity was only brought to the fore in \cite{1007.4330}.  

\paragraph {} P{\'e}rez-Treil-Volberg used the full strength of the non-homogeneous Harmonic analysis, and in particular the innovative paper 
\cite{1003.1596}, to prove a remarkable extension of the $ T1$ Theorem to the $ A_2$ setting.  Loosely, an operator $ T$ with a Calder\'on-Zygmund kernel, then $ T$ extends to a bounded operator on $ L ^{2} (w)$, $ w\in A_2$, if and only if the testing conditions of Theorem~\ref{t.well} hold.    Then, it was shown \cite{1006.2530} that the linear bound holds for Calder\'on-Zygmund operators with sufficiently smooth kernels.  This proof used the Belykin-Coifman-Roklin \cite{MR1085827} decomposition, and the method of \cite{MR2657437} to verify the testing conditions.  
 A short time later, Hyt{\"o}nen \cite{1007.4330}, used a random variant of the  Belykin-Coifman-Roklin method to give a proof of the linear bound for arbitrary smoothness, again using the $ A_2$ $ T1$ Theorem of \cite{1003.1596}. This proof of the full conjecture was then streamlined in \cite{1010.0755}, giving the line of argument we have followed in this survey. 

\paragraph {}  Lerner has conjectured that the weak-type bound on Calder\'on-Zygmund operators should obey the linear bound in $ A_p$ for all $ 1<p< \infty$.  This has been verified for dyadic Calder\'on-Zygmund operators, without careful attention to behavior of the exponents in terms of complexity \cite{0911.0713}, and for the smooth case, with enough derivatives, in \cite{1006.2530}.   The principal technique is again derived from \cite{MR2657437}, as well as a (simple) testing condition for the weak-type inequality for singular integrals given in \cite{0805.0246}, also see \cite{0911.3920}. 
Indeed, this argument proves the linear bound in $ A_p$ for the \emph{maximal truncations} of singular integrals, as this is the kind of operator that we have the testing conditions for.  
It seems likely that this conjecture would follow from Theorem~\ref{t.bcr}, if one tracks complexity constants.    

\paragraph {}  The endpoint case of these estimates is also of interest, namely, for $ p=1$.  It is an elementary consequence of a covering lemma argument that for an arbitrary weight $ w$, the Maximal Function $ M$ maps $ L ^{ 1} (M w )$ into $ L ^{1,\infty } (w)$.  
It was then the subject of conjecture if the same inequality holds for singular integrals.  This was disproved for Haar multipliers by Maria Reguera  \cite{1008.3943}, and then for the Hilbert transform by Reguera-Thiele \cite{1011.1767}.   

\paragraph {} With the failure of the most optimistic form of the conjecture above, one can then ask if its natural variant for $ w\in A_1$ holds.  Namely, does the Hilbert transform map $ L ^{1} (w)$ into $ L ^{1, \infty } (w) $ for $ w\in A_1$, with norm estimate dominated by a constant times $ \lVert w\rVert_{A_1}$?  This also fails in the dyadic case \cite{nrvv}.  On the other hand, the Hilbert transform does 
map  $ L ^{1} (w)$ into $ L ^{1, \infty } (w) $, and the best known upper bound on the norm is $ \lVert w\rVert_{A_1} \log _{+} \lVert w\rVert_{A_1}$.  See \cites{MR2480568} for more information on these last two points.   

\paragraph {} A interesting part of the linear bound in $ A_2$ is that one needs a substantive portion of two-weight theory to address it. 
This is Theorem~\ref{t.well} above.  The general two-weight question is a rather intricate one, with a full discussion carrying us beyond the scope of this text.  The interested reader should consult \cite{V} for a general introduction, and the more recent papers \cites{0911.3920,0805.0246,1001.4043,1003.1596}.


\begin{bibsection}
\begin{biblist}
\bib{MR1085827}{article}{
  author={Beylkin, G.},
  author={Coifman, R.},
  author={Rokhlin, V.},
  title={Fast wavelet transforms and numerical algorithms. I},
  journal={Comm. Pure Appl. Math.},
  volume={44},
  date={1991},
  number={2},
  pages={141--183},
  issn={0010-3640},
}

\bib{MR2433959}{article}{
  author={Beznosova, Oleksandra V.},
  title={Linear bound for the dyadic paraproduct on weighted Lebesgue space $L\sb 2(w)$},
  journal={J. Funct. Anal.},
  volume={255},
  date={2008},
  number={4},
  pages={994--1007},
  issn={0022-1236},
}

\bib{MR1124164}{article}{
  author={Buckley, Stephen M.},
  title={Estimates for operator norms on weighted spaces and reverse Jensen inequalities},
  journal={Trans. Amer. Math. Soc.},
  volume={340},
  date={1993},
  number={1},
  pages={253--272},
  issn={0002-9947},
}

\bib{1002.2396}{article}{
  author={{Chung}, Daewon},
  author={P{\'e}rez, Carlos},
  author={Pereyra, Mar{\'{\i }}a Cristina},
  title={Sharp bounds for general commutators on weighted Lebesgue spaces},
  eprint={http://arxiv.org/abs/1002.2396},
  date={2010},
}

\bib{MR2628851}{article}{
  author={Cruz-Uribe, David},
  author={Martell, Jos{\'e} Mar{\'{\i }}a},
  author={P{\'e}rez, Carlos},
  title={Sharp weighted estimates for approximating dyadic operators},
  journal={Electron. Res. Announc. Math. Sci.},
  volume={17},
  date={2010},
  pages={12--19},
  issn={1935-9179},
}

\bib{1001.4254}{article}{
  author={Cruz-Uribe, David},
  author={Martell, Jos{\'e} Mar{\'{\i }}a},
  author={P{\'e}rez, Carlos},
  title={Sharp weighted estimates for classical operators},
  date={2010},
  eprint={http://arxiv.org/abs/1001.4254},
}

\bib{MR763911}{article}{
  author={David, Guy},
  author={Journ{\'e}, Jean-Lin},
  title={A boundedness criterion for generalized Calder\'on-Zygmund operators},
  journal={Ann. of Math. (2)},
  volume={120},
  date={1984},
  number={2},
  pages={371--397},
  issn={0003-486X},
}

\bib{MR1992955}{article}{
  author={Dragi{\v {c}}evi{\'c}, Oliver},
  author={Volberg, Alexander},
  title={Sharp estimate of the Ahlfors-Beurling operator via averaging martingale transforms},
  journal={Michigan Math. J.},
  volume={51},
  date={2003},
  number={2},
  pages={415--435},
  issn={0026-2285},
  review={\MR {1992955 (2004c:42030)}},
  doi={10.1307/mmj/1060013205},
}

\bib{MR2140200}{article}{
  author={Dragi{\v {c}}evi{\'c}, Oliver},
  author={Grafakos, Loukas},
  author={Pereyra, Mar{\'{\i }}a Cristina},
  author={Petermichl, Stefanie},
  title={Extrapolation and sharp norm estimates for classical operators on weighted Lebesgue spaces},
  journal={Publ. Mat.},
  volume={49},
  date={2005},
  number={1},
  pages={73--91},
  issn={0214-1493},
}

\bib{MR1439553}{article}{
  author={Fefferman, R.},
  author={Pipher, J.},
  title={Multiparameter operators and sharp weighted inequalities},
  journal={Amer. J. Math.},
  volume={119},
  date={1997},
  number={2},
  pages={337--369},
  issn={0002-9327},
}

\bib{MR1110189}{article}{
  author={Figiel, Tadeusz},
  title={Singular integral operators: a martingale approach},
  conference={ title={Geometry of Banach spaces}, address={Strobl}, date={1989}, },
  book={ series={London Math. Soc. Lecture Note Ser.}, volume={158}, publisher={Cambridge Univ. Press}, place={Cambridge}, },
  date={1990},
  pages={95--110},
}

\bib{MR0312139}{article}{
  author={Hunt, Richard},
  author={Muckenhoupt, Benjamin},
  author={Wheeden, Richard},
  title={Weighted norm inequalities for the conjugate function and Hilbert transform},
  journal={Trans. Amer. Math. Soc.},
  volume={176},
  date={1973},
  pages={227--251},
  issn={0002-9947},
}

\bib{MR2464252}{article}{
  author={Hyt{\"o}nen, Tuomas},
  title={On Petermichl's dyadic shift and the Hilbert transform},
  language={English, with English and French summaries},
  journal={C. R. Math. Acad. Sci. Paris},
  volume={346},
  date={2008},
  number={21-22},
  pages={1133--1136},
  issn={1631-073X},
}

\bib{1007.4330}{article}{
  author={Hyt\"onen, Tuomas},
  title={The sharp weighted bound for general Calderon-Zygmund operators},
  eprint={http://arxiv.org/abs/1007.4330},
  date={2010},
}

\bib{0911.0713}{article}{
  author={Hyt\"onen, Tuomas},
  author={Lacey, Michael T.},
  author={Reguera, Maria Carmen},
  author={Vagharshakyan, Armen},
  title={Weak and Strong-type estimates for Haar Shift Operators: Sharp power on the $A_p$ characteristic},
  eprint={http://www.arxiv.org/abs/0911.0713},
  date={2009},
}

\bib{1006.2530}{article}{
  author={Hyt\"onen, Tuomas},
  author={Lacey, Michael T.},
  author={Reguera, Maria Carmen},
  author={Sawyer, Eric T.},
  author={Uriarte-Tuero, Ignacio},
  author={Vagharshakyan, Armen},
  title={Weak and Strong type $ A_p$ Estimates for Calder—n-Zygmund Operators},
  eprint={http://www.arxiv.org/abs/1006.2530},
  date={2010},
}

\bib{1010.0755}{article}{
  author={Hyt\"onen, T.},
  author={P{\'e}rez, Carlos},
  author={Treil, S.},
  author={Volberg, A.},
  title={Sharp weighted estimates of the dyadic shifts and $A_2$ conjecture},
  journal={ArXiv e-prints},
  eprint={http://arxiv.org/abs/1010.0755},
}

\bib{gk}{article}{
  author={Gabidzashvili, M. A.},
  author={Kokilashvili, V.},
  title={Two weight weak type inequalities for fractional type integrals},
  journal={Ceskoslovenska Akademie Ved.},
  volume={45},
  date={1989},
  pages={1--11},
}

\bib{MR1791462}{book}{
  author={Genebashvili, Ioseb},
  author={Gogatishvili, Amiran},
  author={Kokilashvili, Vakhtang},
  author={Krbec, Miroslav},
  title={Weight theory for integral transforms on spaces of homogeneous type},
  series={Pitman Monographs and Surveys in Pure and Applied Mathematics},
  volume={92},
  publisher={Longman},
  place={Harlow},
  date={1998},
  pages={xii+410},
  isbn={0-582-30295-1},
  review={\MR {1791462 (2003b:42002)}},
}

\bib{kwon-2007}{article}{
  author={Kwon, Hyun-Kyoung},
  author={Treil, S.},
  title={Similarity of operators and geometry of eigenvector bundles},
  eprint={http://www.arxiv.org/abs/0712.0114},
  date={2007},
}

\bib{MR2652182}{article}{
  author={Lacey, Michael T.},
  author={Moen, Kabe},
  author={P{\'e}rez, Carlos},
  author={Torres, Rodolfo H.},
  title={Sharp weighted bounds for fractional integral operators},
  journal={J. Funct. Anal.},
  volume={259},
  date={2010},
  number={5},
  pages={1073--1097},
  issn={0022-1236},
}

\bib{MR2657437}{article}{
  author={Lacey, Michael T.},
  author={Petermichl, Stefanie},
  author={Reguera, Maria Carmen},
  title={Sharp $A_2$ inequality for Haar shift operators},
  journal={Math. Ann.},
  volume={348},
  date={2010},
  number={1},
  pages={127--141},
  issn={0025-5831},
}

\bib{0805.4711}{article}{
  author={Lacey, Michael T.},
  author={Sawyer, Eric T.},
  author={Uriarte-Tuero, Ignacio},
  title={Astala's Conjecture on Distortion of Hausdorff Measures under Quasiconformal Maps in the Plane},
  date={2008},
  journal={Acta Math., to appear},
  eprint={http://www.arxiv.org/abs/0805.4711},
}

\bib{0911.3437}{article}{
  author={Lacey, Michael T.},
  author={Sawyer, Eric T.},
  author={Uriarte-Tuero, Ignacio},
  title={Two Weight Inequalities for Discrete Positive Operators},
  date={2009},
  journal={Submitted},
  eprint={http://www.arxiv.org/abs/0911.3437},
}

\bib{0805.0246}{article}{
  author={Lacey, Michael T.},
  author={Sawyer, Eric T.},
  author={Uriarte-Tuero, Ignacio},
  title={A characterization of two weight norm inequalities for maximal singular integrals with one doubling measure},
  date={2008},
  journal={Submitted to Analysis and PDE.},
  eprint={http://arxiv.org/abs/0805.0246},
}

\bib{0911.3920}{article}{
  author={Lacey, Michael T.},
  author={Sawyer, Eric T.},
  author={Uriarte-Tuero, Ignacio},
  title={Two Weight Inequalities for Maximal Truncations of Dyadic Calder\'on-Zygmund Operators},
  date={2009},
  journal={Submitted},
  eprint={http://www.arxiv.org/abs/0911.3920},
}

\bib{1001.4043}{article}{
  author={Lacey, Michael T.},
  author={Sawyer, Eric T.},
  author={Uriarte-Tuero, Ignacio},
  title={A Two Weight Inequality for the Hilbert transform Assuming an Energy Hypothesis},
  eprint={http://www.arXiv.org/abs/1001.4043},
  date={2010},
}

\bib{localSharp}{article}{
  author={Lerner, Andrei K.},
  title={A pointwise estimate for local sharp maximal function with applications to singular integrals},
  journal={Bull. LMS, to appear},
  date={2009},
}

\bib{1005.1422}{article}{
  author={Lerner, Andrei K.},
  title={Sharp weighted norm inequalities for Littlewood-Paley operators and singular integrals},
  date={2010},
  eprint={http://arxiv.org/abs/1005.1422},
}

\bib{MR2399047}{article}{
  author={Lerner, Andrei K.},
  title={An elementary approach to several results on the Hardy-Littlewood maximal operator},
  journal={Proc. Amer. Math. Soc.},
  volume={136},
  date={2008},
  number={8},
  pages={2829--2833},
  issn={0002-9939},
}

\bib{MR2480568}{article}{
  author={Lerner, Andrei K.},
  author={Ombrosi, Sheldy},
  author={P{\'e}rez, Carlos},
  title={$A\sb 1$ bounds for Calder\'on-Zygmund operators related to a problem of Muckenhoupt and Wheeden},
  journal={Math. Res. Lett.},
  volume={16},
  date={2009},
  number={1},
  pages={149--156},
  issn={1073-2780},
}

\bib{MR0293384}{article}{
  author={Muckenhoupt, Benjamin},
  title={Weighted norm inequalities for the Hardy maximal function},
  journal={Trans. Amer. Math. Soc.},
  volume={165},
  date={1972},
  pages={207--226},
  issn={0002-9947},
}

\bib{MR0311856}{article}{
  author={Muckenhoupt, Benjamin},
  title={Hardy's inequality with weights},
  note={Collection of articles honoring the completion by Antoni Zygmund of 50 years of scientific activity, I},
  journal={Studia Math.},
  volume={44},
  date={1972},
  pages={31--38},
  issn={0039-3223},
  review={\MR {0311856 (47 \#418)}},
}

\bib{nrvv}{article}{
  author={Nazarov, F.},
  author={Reznikov, Alexander},
  author={Vasyunin, Vasily},
  author={Volberg, Alexander},
  title={Personal Communication},
  date={2010},
}

\bib{MR1998349}{article}{
  author={Nazarov, F.},
  author={Treil, S.},
  author={Volberg, A.},
  title={The $Tb$-theorem on non-homogeneous spaces},
  journal={Acta Math.},
  volume={190},
  date={2003},
  number={2},
  pages={151--239},
}

\bib{1003.1596}{article}{
  author={Nazarov, F.},
  author={Treil, S.},
  author={Volberg, A.},
  title={ Two weight estimate for the Hilbert transform and Corona decomposition for non-doubling measures},
  date={2004},
  eprint={http://arxiv.org/abs/1003.1596},
}

\bib{MR2407233}{article}{
  author={Nazarov, F.},
  author={Treil, S.},
  author={Volberg, A.},
  title={Two weight inequalities for individual Haar multipliers and other well localized operators},
  journal={Math. Res. Lett.},
  volume={15},
  date={2008},
  number={3},
  pages={583--597},
  issn={1073-2780},
}

\bib{MR1685781}{article}{
  author={Nazarov, F.},
  author={Treil, S.},
  author={Volberg, A.},
  title={The Bellman functions and two-weight inequalities for Haar multipliers},
  journal={J. Amer. Math. Soc.},
  volume={12},
  date={1999},
  number={4},
  pages={909--928},
  issn={0894-0347},
  review={\MR {1685781 (2000k:42009)}},
}

\bib{NV}{article}{
  author={Nazarov, F.},
  author={Volberg, A.},
  title={The Bellman function, the two-weight Hilbert transform, and embeddings of the model spaces $K_\theta $},
  note={Dedicated to the memory of Thomas H.\ Wolff},
  journal={J. Anal. Math.},
  volume={87},
  date={2002},
  pages={385--414},
}

\bib{MR2465289}{article}{
  author={Nazarov, F.},
  author={Peherstorfer, F.},
  author={Volberg, A.},
  author={Yuditskii, P.},
  title={Asymptotics of the best polynomial approximation of $\vert x\vert \sp p$ and of the best Laurent polynomial approximation of ${\rm sgn}(x)$ on two symmetric intervals},
  journal={Constr. Approx.},
  volume={29},
  date={2009},
  number={1},
  pages={23--39},
  issn={0176-4276},
}

\bib{MR1756958}{article}{
  author={Petermichl, Stefanie},
  title={Dyadic shifts and a logarithmic estimate for Hankel operators with matrix symbol},
  language={English, with English and French summaries},
  journal={C. R. Acad. Sci. Paris S\'er. I Math.},
  volume={330},
  date={2000},
  number={6},
  pages={455--460},
  issn={0764-4442},
}

\bib{MR1964822}{article}{
  author={Petermichl, S.},
  author={Treil, S.},
  author={Volberg, A.},
  title={Why the Riesz transforms are averages of the dyadic shifts?},
  booktitle={Proceedings of the 6th International Conference on Harmonic Analysis and Partial Differential Equations (El Escorial, 2000)},
  journal={Publ. Mat.},
  date={2002},
  number={Vol. Extra},
  pages={209--228},
  issn={0214-1493},
  review={\MR {1964822 (2003m:42028)}},
}

\bib{MR2354322}{article}{
  author={Petermichl, Stefanie},
  title={The sharp bound for the Hilbert transform on weighted Lebesgue spaces in terms of the classical $A\sb p$ characteristic},
  journal={Amer. J. Math.},
  volume={129},
  date={2007},
  number={5},
  pages={1355--1375},
  issn={0002-9327},
}

\bib{MR2367098}{article}{
  author={Petermichl, Stefanie},
  title={The sharp weighted bound for the Riesz transforms},
  journal={Proc. Amer. Math. Soc.},
  volume={136},
  date={2008},
  number={4},
  pages={1237--1249},
  issn={0002-9939},
}

\bib{MR1894362}{article}{
  author={Petermichl, Stefanie},
  author={Volberg, Alexander},
  title={Heating of the Ahlfors-Beurling operator: weakly quasiregular maps on the plane are quasiregular},
  journal={Duke Math. J.},
  volume={112},
  date={2002},
  number={2},
  pages={281--305},
  issn={0012-7094},
}

\bib{1008.3943}{article}{
  author={Reguera, Maria Carmen},
  title={On Muckenhoupt-Wheeden Conjecture},
  eprint={http://www.arxiv.org/abs/1008.3943},
  date={2010},
}

\bib{1011.1767}{article}{
  author={Reguera, Maria Carmen},
  author={Thiele, Christoph},
  title={The Hilbert transform does not map $L^1(Mw)$ to $L^{1,\infty }(w)$},
  eprint={http://www.arxiv.org/abs/1011.1767},
  date={2010},
}

\bib{MR745140}{article}{
  author={Rubio de Francia, Jos{\'e} L.},
  title={Factorization theory and $A_{p}$ weights},
  journal={Amer. J. Math.},
  volume={106},
  date={1984},
  number={3},
  pages={533--547},
  issn={0002-9327},
}

\bib{MR676801}{article}{
  author={Sawyer, Eric T.},
  title={A characterization of a two-weight norm inequality for maximal operators},
  journal={Studia Math.},
  volume={75},
  date={1982},
  number={1},
  pages={1--11},
  issn={0039-3223},
}

\bib{MR930072}{article}{
  author={Sawyer, Eric T.},
  title={A characterization of two weight norm inequalities for fractional and Poisson integrals},
  journal={Trans. Amer. Math. Soc.},
  volume={308},
  date={1988},
  number={2},
  pages={533--545},
  issn={0002-9947},
}

\bib{MR1175693}{article}{
  author={Sawyer, E.},
  author={Wheeden, R. L.},
  title={Weighted inequalities for fractional integrals on Euclidean and homogeneous spaces},
  journal={Amer. J. Math.},
  volume={114},
  date={1992},
  number={4},
  pages={813--874},
  issn={0002-9327},
  review={\MR {1175693 (94i:42024)}},
  doi={10.2307/2374799},
}

\bib{1007.3089}{article}{
  author={{Scurry}, J.},
  title={A Characterization of Two-Weight Inequalities for a Vector-Valued Operator},
  eprint={http://arxiv.org/abs/1007.3089},
  date={2010},
}

\bib{MR1232192}{book}{
  author={Stein, Elias M.},
  title={Harmonic analysis: real-variable methods, orthogonality, and oscillatory integrals},
  series={Princeton Mathematical Series},
  volume={43},
  note={With the assistance of Timothy S. Murphy; Monographs in Harmonic Analysis, III},
  publisher={Princeton University Press},
  place={Princeton, NJ},
  date={1993},
  pages={xiv+695},
  isbn={0-691-03216-5},
}

\bib{0911.4968}{article}{
  author={{Vagharshakyan}, Armen},
  title={Recovering Singular Integrals from Haar Shifts},
  journal={ArXiv e-prints},
  eprint={http://arxiv.org/abs/1007.2994},
  date={2009},
}

\bib{V}{book}{
  author={Volberg, A.},
  title={Calder\'on-Zygmund capacities and operators on nonhomogeneous spaces},
  series={CBMS Regional Conference Series in Mathematics},
  volume={100},
  publisher={Published for the Conference Board of the Mathematical Sciences, Washington, DC},
  date={2003},
  pages={iv+167},
  isbn={0-8218-3252-2},
}

\bib{MR1748283}{article}{
  author={Wittwer, Janine},
  title={A sharp estimate on the norm of the martingale transform},
  journal={Math. Res. Lett.},
  volume={7},
  date={2000},
  number={1},
  pages={1--12},
  issn={1073-2780},
  review={\MR {1748283 (2001e:42022)}},
}

\end{biblist}
\end{bibsection}
\end{document}